\documentclass[a4paper]{article}

\usepackage{amsmath}
\usepackage{amssymb}
\usepackage{amsfonts}

\usepackage{vmargin}

\newcommand{\tens}{\otimes}
\newcommand{\mat}{\mathcal}

\newcommand{\K}{\ensuremath{\mathbb{K}}}

\newcommand{\B}{\ensuremath{\mathbb{B}}}
\newcommand{\R}{\ensuremath{\mathbb{R}}}
\newcommand{\C}{\ensuremath{\mathbb{C}}}

\newcommand{\Id}{\ensuremath{\mathrm{Id}}}
\newcommand{\vect}{\ensuremath{\mathop{\rm Span\,}\nolimits}}
\newcommand{\im}{\ensuremath{\mathop{\rm Ran\,}}}
\renewcommand{\dim}{\ensuremath{\mathop{\rm dim\,}}}

\newcommand{\mk}{\medskip}

\renewcommand{\le}{\ensuremath{\leqslant}}
\renewcommand{\ge}{\ensuremath{\geqslant}}
\renewcommand{\leq}{\ensuremath{\leqslant}}
\renewcommand{\geq}{\ensuremath{\geqslant}}

\renewcommand{\u}{\underline}
\renewcommand{\o}{\overline}
\renewcommand{\b}{\mathring}
\newcommand{\la}{\langle}
\newcommand{\ra}{\rangle}
\font\bbold=bbold12
\newcommand{\ind}{\textrm{\bbold 1}}
\newcommand{\qed}{\hfill \vrule height6pt  width6pt depth0pt}
\newtheorem{thm}{Theorem}[section]

\newtheorem{prop}[thm]{Proposition}
\newtheorem{cor}[thm]{Corollary}
\newtheorem{lemma}[thm]{Lemma}

\newenvironment{rk}{\refstepcounter{thm}\noindent%
{\bf Remark \arabic{section}.\arabic{thm}} \ }{
\smallskip

}

\newenvironment{pf}[1][]{\noindent {\it Proof #1} : }{\hbox{~}\qed
\smallskip
}

\title{The von Neumann algebra generated by  $t$-gaussians}
\date{}
\author{\'Eric Ricard}

\begin{document}

\maketitle

\begin{abstract}
 We study the $t$-deformation of gaussian von Neumann algebras. 
When the number of generators is fixed, it is proved
that if $t$ sufficiently close to 1, then these algebras do not depend on $t$.
In the same way, the notion of conditionally free von Neumann algebras often 
coincides with freeness.
\end{abstract}

\makeatletter
\renewcommand{\@makefntext}[1]{#1}
\makeatother
\footnotetext{\noindent Laboratoire de Math\'{e}matiques,
 Universit\'{e} de Franche-Comt\'{e},
 25030 Besan\c con, cedex - France \\ eric.ricard@math.univ-fcomte.fr}

\section{Introduction}

~
 
 \indent The notions of free and non commutative probabilities originally
appeared  in the works of 
Voiculescu in the 80's (see \cite{VDN} for instance) to study von 
Neumann algebras, in particular the von Neumann algebra $L(\mathbb F_n)$
associated to the free group with $n$ generators. Since then, this
domain has expanded rapidly, and is now considered as a subject in itself.
It has connexions with combinatorics, classical probabilities, mathematical 
physics and of course operator algebras. Naturally, people become 
interested in finding  
generalizations or deformations of the free probability or the free product
constructions to obtain new non commutative probability spaces,
especially from the combinatoric point of view. One of the first successful
attempt was made by Bo\.zejko and Speicher (\cite{BS2}, see also \cite{BKS})
in introducing the so called
$q$-deformation of the free factor. These von Neumann algebras were 
studied in the last years, at the moment it is known that they share some
properties with the free group algebras : they are factor \cite{R},
 non injective \cite{N}, 
and solid \cite{Sh} 
(for some values of the parameter $q$). Nevertheless, it is still
unknown if they really differ from the free group algebras.

With the same spirit, 
Bo\.zejko, Leinert and Speicher introduced in \cite{BLS} the concept 
of conditional freeness. They were able to describe the combinatoric 
underlying this notion in a way similar to that of the freeness in terms
of Cauchy Transforms and their reciprocal.  
The school of 
Accardi also developed its own deformation by considering the algebra generated
by position operators on an interacting Fock space
(see \cite{AB}). Unfortunately, very 
little work on these two topics were concerned with the study of the 
resulting von Neumann object. One of the motivation for the present work is 
to start such a study like in the $q$-case.
 We will mainly focus on specific examples of  von 
Neumann algebras, the $t$-gaussian von Neumann algebras. They have the 
advantage to be both a model for conditional freeness and an algebra on 
an interacting Fock space. 
They even appear as a limit object for the theory of 
the $t$-convolution of measures of Bo\.zejko and Wysocza{\'n}ski
in \cite{BW}. So they seem to be quite central in all those deformations of 
the free probability.        

 As an  illustration, we recall the definition of the conditional freeness,
and explain how the $t$-gaussian algebra appears. Let
$(A_i,\phi_i,\psi_i)$ be $*$-algebras equipped with a pair of states.
Assume that all $A_i$'s lie in a bigger algebra $A$ with a state
$\phi$. Then we say that the algebras $A_i$'s are conditionally free for 
$\phi$ provided that whenever 
$$a_j\in A_{i_j},\; i_1\neq i_2\neq ...\neq i_n, \; \psi_{i_j}(a_j)=0$$
we have
$$\phi(a_{i_1}...a_{i_n})=\phi_{i_1}(a_{i_1})...\phi_{i_n}(a_{i_n})$$
 It is clear that for any given  $(A_i,\phi_i,\psi_i)$  it is possible to 
construct $\phi$ on the free product $A=*_{i=1}^nA_i$ so that the
$A_i$'s are conditionally free. In the situation where $\phi_i=\psi_i$, one
recovers the classical notion of freeness, that is $\phi$ is the free
product of the $\psi_i$. We let $\psi$ be the free product of the 
$\psi_i$ on $A$.
 
Given  some probability measures $\mu_i,\,\nu_i$ (with compact support), 
one can consider them as states on the space of polynomials in one variable 
$\C[X_i]$ in a natural manner by the formula
$$\mu_i(X_i^k)=\int x^k {\rm d} \mu_i(x).$$ 
 The free product of $n$ algebras of polynomials  in one variable is exactly 
the space of polynomials in $n$ non commuting variables $\C\la X_i\ra_{i\le n}$. So, using the above construction, from the measures  $\mu_i,\,\nu_i$ one can
build two new states $\phi$ and $\psi$ on  $\C\la X_i\ra$, 
so that the algebras   
$(\C[ X_i],\mu_i,\nu_i)$ are conditionally free in 
$(\C\la X_i\ra,\phi,\psi)$.
 The  distribution of $X_1+X_2$ with respect to $\phi$ is the compactly 
supported measure $\mu$
so that
$$\phi((X_1+X_2)^k)=\int x^k {\rm d}\mu(x).$$  
 The c-free convolution of 
$(\mu_1,\nu_1)$ and  $(\mu_2,\nu_2)$ is then defined as 
$$(\mu_1,\nu_1)\boxplus_c (\mu_2,\nu_2) =(\mu, \nu_1 \boxplus  \nu_2)$$
where  $\boxplus$ means the usual free convolution (i.e. $ \nu_1 \boxplus  
\nu_2$ is 
the distribution of 
 $X_1+X_2$ with respect to $\psi$). This operation
$\boxplus_c $ on couples of measures is commutative and associative.
Limits theorems were obtained in this framework. 
 For instance, a central limit theorem consists in finding the possible limit 
distributions of $X_1+...+X_n/\sqrt{n}$, where the $X_i$ are identically 
distributed (and centered) and are conditionally free. 
The measures $(\mu,\nu)$ appearing at the limit, 
are parameterised by the second moments
of $\mu$ and $\nu$. So up to normalization $\mu_t(X^2)=1$ and $\nu_t(X^2)=t$.
 The von Neumann algebra
arising from the GNS 
construction of the first state on the conditional free product 
$*_{i=1}^n(\C[ X_i],\mu_t,\nu_t)$, is $\Gamma_{t,n}$ the $t$-gaussian
algebra with $n$ generators. When $t=1$, this is just the classical objects 
of free probabilities, that is $\Gamma_{1,n}=L(\mathbb F_n)$. 
 
  The main result of this paper states that for $n\ge 2$ 
$$\Gamma_{t,n}=\left\{
\begin{array}{cl}
\Gamma_{1,n} & \textrm{ if }\;t\in
\Big[\frac n {n+\sqrt n},\frac n {n-\sqrt n}\Big]\\
\Gamma_{1,n} \oplus \B(\ell_2) & \textrm{ otherwise }
\end{array}\right.,$$
where as usual $\B(\ell_2)$ stands for the bounded operators on a separable
Hilbert space. Consequently, this result is a bit disappointing for 
operator algebras, as the deformations may not produce new objects.

 In the next section, we give a precise construction of $\Gamma_{t,n}$ arising 
from the theory of one mode interacting Fock space and very basic results about
it. The third section is devoted to the proof of the main result. 
 In the last section, we give some conditions that ensure the equality between 
the reduced free product of commutative von Neumann algebras $(A_i,\psi_i)$
 and their conditional reduced free product $(A_i,\phi_i,\psi_i)$. These 
conditions are satisfied for the $t$-gaussian algebras (in the first case),
however, we prefer to give slightly different proofs for these examples as one 
can extract information from them to get some partial 
results on interacting Fock space.

\section{Basics}
 
 In the whole paper, $t$ will be a positive real number $t>0$. We will also use
standard notations, $\B(H)$ and $\K(H)$ will denote the bounded and the compact
operators on the Hilbert space $H$.

 For a given Hilbert space $\mat H$, with real part
 $\mat H_\R$, equipped with a scalar product $\la .,.\ra$, we denote by
$\mat F_1$, the Fock space built on $\mat H$ :
$$\mat F_1= \C \Omega\oplus_{k\ge1}  \mat H^{\tens k}$$
 More generally, the $t$-deformed Fock space is given by 
$$\mat F_t=\C \Omega\oplus_{k\ge1} t^{k-1} \mat H^{\tens k}$$
Here, $ t^{k-1} \mat H^{\tens k}$ simply means the space 
$H^{\tens k}$ where the scalar product is multiplied by $t^{k-1}$.
  
 In the most part of the paper, we will assume that $\mat H$ is finite 
dimensional, say $\dim \mat H=n$. Let $(e_i)_{i=1...n}$  
be a real basis of $\mat H_\R$. Then one can define
a canonical basis for $\mat F_t$. Its elements $e_{\u i}$ are indexed by 
words $\u i$ in the $n$ letters $1,...,n$, and
$$e_{\u i}=\frac 1 {{\sqrt t}^{k-1} }
e_{i_1}\tens...\tens e_{i_k} \qquad \textrm{ for  }\, \u i =i_1...i_k$$ 

 As usual, for $e\in \mat H_\R$, the creation operator associated to $e$ is 
defined on $\mat F_t$ by 
$$\begin{array}{l}
l_t(e)\Omega=e \\
l_t(e)(h_1\tens...\tens h_k)=e\tens h_1\tens...\tens h_k
\end{array}$$
It is well known that $l_t(e)$ extends to a bounded operator on $\mat F_t$.
Its adjoint is given by 
$$\begin{array}{l}
l_t(e)^*\Omega=0 \\
l_t(e)^*h=\la h,e\ra \Omega \\ 
l_t(e)^*(h_1\tens...\tens h_k)=t <h_1,e>\, h_2\tens...\tens h_k
\end{array}$$
The $t$-gaussian associated to $e$ is the operator $s^t(e)=l_t(e)+l_t(e)^*$.
 
 We are interested in the von Neumann generated by the $t$-gaussians.
Let $s_i^t=l_t(e_i)+l_t(e_i)^*$ for $i=1,...n$.
$C_{t,n}$ and $\Gamma_{t,n}$ refer to the $C^*$-algebra and the 
von Neumann algebras
generated by the $s_i^t$'s, $i=1,...,n$ :
$$\Gamma_{t,n}=\{ l_t(e)+l_t(e)^* \,;\, e\in \mat H_\R \} ''\subset \B(\mat 
F_t)$$

 This deformation of the Fock Hilbert space is a particular case of  
interactive Fock spaces. This type of objects was introduced in \cite{AB}.
The basic idea is to modify the scalar product at each level 
$\mat H^{\tens k}$ by a positive scalar $\lambda_k$. 

 The vector state $\la . \Omega,\Omega \ra$ on $\B(\mat F_t)$ is called the
vacuum state and will be denoted by $\phi$. In the case $t=1$, $\phi$ is a
trace on $\Gamma_{1,n}$, so we will prefer the notation $\tau$.

\medskip

 In the following, we sum up all basic results about the von Neumann algebra
generated by a single $t$-gaussian $s^t$ when $\mat H=\C$. 
\begin{prop}\label{basic}
 $\Gamma_{t,1}$ is in GNS position with respect to $\phi$, which is faithful
on it.\\
 The distribution of $s^t$ with respect to $\phi$ is given by 
$$\begin{array}{ll}
\displaystyle{
\frac 1 {2\pi} \frac{\sqrt{4t-x^2}}{1-(1-t)x^2} \ind_{[-2\sqrt t,2\sqrt t]} 
} {\rm d}x & \textrm {if \; \; $t\geq \frac 1 2$} \\
\displaystyle{
\frac 1 {2\pi} \frac{\sqrt{4t-x^2}}{1-(1-t)x^2} \ind_{[-2\sqrt t,2\sqrt t]}
{\rm d}x
+\,\frac{1-2t}{2-2t} (\delta_{\frac 1 {-\sqrt{1-t}}}+
\delta_{\frac 1 {\sqrt{1-t}}})} & \textrm {if \; \;  $t< \frac 1 2$}
\end{array}$$
The map $\rho:\Gamma_{t,1}\to \Gamma_{1,1}$ given by $\rho(s^t)=\sqrt t s^1$
 extends 
to a normal representation. \\
Moreover, given $i\leq n$, 
there are $*$-isometric normal representations $\pi_i:\Gamma_{t,1}
\to \Gamma_{t,n}$ given by $\pi(s^t)=s_i^t$.

\end{prop} 

\begin{pf}
 The computations of the distribution of $s^t$ can be found in 
\cite{BLS,BW,Wo}. The 
$G$-transform of the distribution of $s^t$ for $\phi$ is given by 
$$G_{s^t}(z)=\frac {(\frac 1 2 -t)z+\frac 1 2 \sqrt{z^2-4t}}{z^2(1-t)-1}
=\cfrac 1 {z - \cfrac 1 {z - \cfrac t {z -  \cfrac t {z -\cfrac
 t \ddots}}}}$$  
It is obvious that $\Omega$ is cyclic for $\Gamma_{t,1}$ so is also separating
since the algebra is commutative and $\phi$ is faithful.
 In particular, the spectral measure of $\sqrt t s^1$ 
is absolutely continuous with respect 
to that of $s^t$,  this implies that the map 
$\rho:\Gamma_{t,1}\to \Gamma_{1,1}$ is well defined and normal.

 Note that the matrix of $s^t$ in the natural orthonormal basis is given by
$$ \left[\begin{array}{ccccccc}
 0 & 1 &  & &  & &  \\
 1 & 0 & \sqrt t & &  &  & \\
   & \sqrt t & 0 &\sqrt t &  & & \\
   &  &\ddots &\ddots &\ddots & & \\
   & &  &\sqrt t & 0 &\sqrt t  &   \\
   & &  &        & \ddots & \ddots &\ddots  \\
   \end{array}\right]$$

 So at the $C^*$-level, $\rho$ is just the restriction of the Calkin map 
($\rho : \B(\mat F_1)\to \B(\mat F_1)/\K(\mat F_1)$, where $\K$ stands for the 
compact operators) since 
$s^t -\sqrt t s^1$ is of rank 2 and $\Gamma_{1,1}$ does not meet the compact
operators (no atoms in the measure).

 In $\Gamma_{t,n}$, $s_i^t$ is unitarily equivalent to $s^t \oplus
(\sqrt t s^1)_\infty$ corresponding to its decomposition in reducing
subspaces (they are indexed by reduced word starting by a letter which
is not a $i$).
\end{pf}

\begin{rk}
 There are natural isometries between all $\mat F_t$ by identifying the 
canonical basis. This allows to consider all $\Gamma_{t,n}$ as subalgebras
of $\B(\mat F_1)$. When we will talk about the Calkin map, we mean the
classical quotient map $\rho : \B(\mat F_1)\to   \B(\mat F_1)/ \K(\mat F_1)$.
Using this identification, we have $\rho(C_{t,n})=\rho(C_{1,n})$ for any 
$0<t$ and $1\leq n <\infty$.
\end{rk}

\begin{rk}
 From the densities, it is clear that as von Neumann algebras, we have
$$\Gamma_{t,1}=\left\{\begin{array}{cl}
 \Gamma_{1,1}\, \oplus\, \C^2 & \textrm{ if } \, t<\frac 12 \\[2pt]
 \Gamma_{1,1}& \textrm{ if } \, t\geq \frac 12 
\end{array}\right.$$
\end{rk}

Let $A_i\approx \Gamma_{t,1}$ be the algebra 
generated by $s^t_i$ in $\Gamma_{i,n}$. By the 
previous proposition,
 there are two given states on $A_i$. One is  coming from the vacuum
state (denoted by $\phi$, no confusion since it coincides with the 
vacuum on $\Gamma_{t,n}$) and another one coming from the vacuum state 
on $\Gamma_{1,1}$ that
is $\tau\rho\pi_i^{-1}$ (call it $\psi$ even if it depends on $i$).

 The following is well known 

\begin{prop}\label{condfree}
The algebras $(A_i,\phi,\psi)$ 
are conditionally free with respect to the vacuum state, that is 
$$\phi(a_1...a_p)=\phi(a_1)...\phi(a_p) \quad \textrm{
whenever $a_j\in A_{i_j}$, 
$i_1\neq i_2\neq...\neq i_p$ and $\psi(a_j)=0$}.$$
Moreover $(\Gamma_{t,n},\phi,\Omega)$ is in GNS position.   
\end{prop}

 We postpone the proof of this fact to the next section.

 Let $U$ be an orthogonal transformation on $\mat H_\R$ and still denote 
its tensorization with $\Id_\C$ by $U$. The first quantization $\Gamma(U)$ 
of $U$ is the unitary on $\mat F_t$ given by 
$$ \Gamma(U)= \Id_{\C \Omega} \oplus_{k\ge 1} U^{\tens^k}$$
Then, $\Gamma_{t,n}$ is stable by conjugation by $\Gamma(U)$ and for any 
$e\in \mat H_\R$,
$$ \Gamma(U) s^t(e)  \Gamma(U)^* =s^t (Ue)$$

\section{The von Neumann algebra $\Gamma_{t,n}$}

\subsection{Factoriality}

 The following is well known to specialists :
\begin{prop}\label{compb}
Let $M\subset \B(H)$ be a von Neumann algebra, if $M$ contains a non
 zero compact operator, then either  $M=\B(K)\tens \Id_d$ with $H=K^d$ and $d$ 
is the smallest rank of a non zero compact projection in $M$ or $M$ is not a 
factor and has a direct summand isomorphic to $\B(K)$ for some $K$.
\end{prop}

\begin{pf}
We know that $M$ contains a finite rank projection. So let $p$ 
be a finite rank projection of minimum rank $d\ge 1$. 
Let $(p_i)_{i\in I}$ be a maximal 
family of mutually orthogonal projections equivalent to $p$ 
(so all $p_i$'s have the same finite rank). 
Let $q=\sum p_i$, we show that $q$ is central.

If $q$ is not central then  there is some $x\in M$ so that $(1-q)xq\neq 0$. 
Hence we can assume that $\im x$ is orthogonal 
to $\im q$ and $xq\neq 0$. There must be an $i$ so that $xp_i\neq 0$.
Consider $xp_i x^*$, it is non zero and of rank less or equal to 
 $d$. By definition 
of $d$, $p_0=xp_ix^*$ must have at most one non zero eigenvalue, and has 
at least one as $p_0\neq 0$. So $p_0$ is a multiple (say 1) of a  projection
of rank $d$. Moreover, $p_ix^*xp_i$ is also a non zero 
projection hence it must be equal to 
$p_i$. So $p_0$ is equivalent to $p$. But the range of $p_0$ is orthogonal 
to $q$ so $p_0q=qp_0=0$, this contradicts the maximality of $I$.     
So $q$ is central.

Let $u_{i,j}$ be the partial isometries between $p_i$ and $p_j$
$i,j\in I$. It is not hard to see that $u_{i,j}$ is a system of matrix units
that generates $qM$. So if $q=1$ then $M$ is a factor, else $M$ is not a 
factor and has a direct summand isomorphic to $\B(\ell_2(I))$.  
\end{pf}  

 In the case $d=1$, the argument is much simpler : let $\xi$ be in
the range of $p$ as above. Then $M$ contains all projections on $K=\overline
{M.\xi}$ (using conjugation of $p$). So if $K\neq H$, then the projection $q$
onto $K$ belongs to both $M'$ (by definition) and $M$ (because $q=\sum p_i$
with $p_i$ projection onto lines corresponding to an onb in $K$), so $M$ is not
a factor and has a direct summand isomorphic to $\B(K)$.
 If $K=H$ the same kind of arguments gives that $M=\B(K)$. 

\subsection{Orthogonal polynomials}

 Let $U_n$ be the Tchebychef polynomial of the second kind of degree $n$.

\begin{prop}
 The orthonormal polynomials for $s^t$ with respect to $\phi$ are given by :
$$ v_0(X)=1  \qquad v_1(X)=X=\sqrt t U_1\Big(\frac X{2\sqrt t}\Big)$$
$$v_n(X)=\sqrt t\left(U_n\Big(\frac X{2\sqrt t}\Big)- \Big(\frac 1 t -1\Big)
U_{n-2}\Big(\frac X{2\sqrt t}\Big)\right) \qquad \textrm{for } n\ge 2$$
The orthonormal polynomials for $\sqrt t s^1$ are given by 
$u_n(X)=U_n\Big(\frac X{2\sqrt t}\Big)$.
\end{prop} 

\begin{pf}
According to the continued fraction 
decomposition of the $G$-transform of the measure, the 
unital orthogonal polynomials must
satisfy the relations 
 $$X P_n(X)= P_{n+1}(X) +t P_{n-1} (X) \qquad \textrm{ for } n\ge 2$$
$$ X P_1(X)= P_2(X) +P_0(X), \qquad P_0(X)=1,\quad P_1(X)=X$$
Then one deduces that this relation is satisfied with 
$P_n(X)=2^n {\sqrt t}^{n-1} v_n(X)$. See also \cite{BF,Wo} for detailed proofs.
\end{pf}

The polynomials $u_n$ and $v_n$ are related one to each other, 
letting $\alpha=\frac 1 t -1$ 
$$
(R)\qquad \qquad \begin{array}{c}
\displaystyle{v_n=\sqrt t \Big( u_n- \alpha u_{n-2}\Big)}\\[8pt]
\displaystyle{u_{2n}=\alpha^n v_0+ \frac 1 
{\sqrt t} \sum_{k=1}^n \alpha^{n-k} v_{2k}}\\[8pt]
\displaystyle{u_{2n+1}=\frac 1 {\sqrt t} \sum_{k=0}^n \alpha^{n-k} v_{2k+1}}
\end{array}$$

\begin{lemma}\label{ident}
At the algebraic level, we have for $i_1\neq...\neq i_l$ 
and $\alpha_j\ge 1$, with $\u i=i_1^{\alpha_1}...i_l^{\alpha_l}$ :
$$ u_{\alpha_1}(s^t_{i_1})...u_{\alpha_{l-1}}(s^t_{i_{l-1}})v_{\alpha_l}(
s^t_{i_l})\Omega= e_{\u i} $$
\end{lemma}

\begin{pf}
 We have $v_{\alpha_l}(s^t_{i_l})\Omega=e_{i_l^{\alpha_l}}$ as the $v_j$ are 
the orthonormal polynomials for $s^t_{i_l}$. Now, $s^t_{i_{l-1}}$ acts on 
tensors starting with a letter $i_l$ exactly as $\sqrt ts^1_{i_{l-1}}$ so 
$u_{\alpha_{l-1}}(s^t_{i_{l-1}})v_{\alpha_l}(
s^t_{i_l})\Omega= e_{i_{l-1}^{\alpha_{l-1}}i_l^{\alpha_l}}$, and so on.
\end{pf}

\begin{pf}[of Proposition \ref{condfree}]
Using multilinearity, we need to show that for 
$i_1\neq...\neq i_l$ 
and $\alpha_j\ge 1$, we have 
$$\phi\big(
u_{\alpha_1}(s^t_{i_1})...u_{\alpha_l}(
s^t_{i_l})\big)=\phi(u_{\alpha_1}(s^t_{i_1}))...\phi(u_{\alpha_l}(
s^t_{i_l}))$$
 From the relations $(R)$, we have 
$$\phi(u_{\alpha_k}(s^t_{i_k}))=
\left\{\begin{array}{ll}\alpha^{\alpha_k/2} & \textrm{if $\alpha_k$ is even}
\\ 0 & \textrm {otherwise }\end{array}\right.$$
But by the previous lemma
$$0=\frac 1 {\sqrt t}\phi\big(u_{\alpha_{1}}(s^t_{i_{1}})...v_{\alpha_l}(
s^t_{i_l})\big)=\phi\big(u_{\alpha_{1}}(s^t_{i_{1}})...u_{\alpha_l}(
s^t_{i_l})\big)-\alpha \phi\big(u_{\alpha_{1}}(s^t_{i_{1}})...u_{\alpha_l-2}(
s^t_{i_l})\big) 
$$
 Then we conclude by induction on $\sum \alpha_k$.

 The statement about the GNS position can then be easily deduced from the 
conditional freeness (see Lemma \ref{baseorth}).
 
\end{pf}

Let $$\eta=\Omega - \sqrt t \alpha \sum_{k=1}^n e_{kk}.$$
Define a functional on $\B(\mat F_t)$ by :
$$\tilde \psi (x)=\la x\Omega, \eta\ra$$
\begin{lemma}\label{densite}
The functional $\tilde \psi$ coincides with $\tau \rho$ on $C_{t,n}$.
\end{lemma}

\begin{pf}
The set of polynomials in the letters $s_i^t$ is dense is 
$C_{t,n}$. As a consequence,  the linear span of products 
$u_{d_1}(s^t_{i_1})...u_{d_m}(s^t_{i_m})$, with $d_j\ge 1$, $m\ge0$
and $i_1\neq i_2\neq... \neq i_m$ is dense in $C_{t,n}$. So we only need
to prove that as soon as $m\ge 1$ (using freeness) :
\begin{eqnarray*}
0=\la u_{d_1}(s^t_{i_1})...u_{d_m}(s^t_{i_m}) \Omega, \eta\ra
&=& \la u_{d_1}(s^t_{i_1})...u_{d_m}(s^t_{i_m})\Omega, 
\Big(1-
\sqrt t \alpha  \sum_{k=1}^n v_2(s^t_k)\Big)\Omega\ra\\
 &=& \phi\left(
\Big(1-
\sqrt t \alpha  \sum_{k=1}^n v_2(s^t_k)\Big)
u_{d_1}(s^t_{i_1})...u_{d_m}(s^t_{i_m})\right)\\
\end{eqnarray*}
So we need to evaluate expressions of the form 
$$\phi( a b_1.... b_l )$$
with $a\in A_{j_0}$, $b_u\in A_{j_{u}}$  with 
$j_0\neq j_1 \neq ... \neq j_{l+1}$ and $\psi(b_i)=0$ (with possibly $l=0$).
Using conditional freeness, and denoting by $b=b_1...b_l$ and 
$x=\phi(b)=\phi(b_1)...\phi(b_l)$

\begin{eqnarray*}
 \phi( a b)&=&\phi \big( (a -\psi(a))b\big) + \psi(a)\phi(b)
\\ & =& \phi(a)\phi(b)
\end{eqnarray*}   
 Using this equality, we get as $\phi(v_2(s^t_i))=0$
\begin{eqnarray*}
\la u_{d_1}(s^t_{i_1})...u_{d_m}(s^t_{i_m}) \Omega, \eta\ra
&=& \phi\left(\Big(1-
\sqrt t \alpha  v_2(s^t_{i_1})\Big)
u_{d_1}(s^t_{i_1})...u_{d_m}(s^t_{i_m})\right)
\end{eqnarray*} 
 In this scalar 
product, there will be a factor 
$$\phi\Big(\Big(1-
\sqrt t \alpha  v_2(s^t_{i_1})\Big)u_{d_1}(s^t_{i_1})\Big)$$
From the formula above, if $d_1$ is odd then $u_{d_1}$ is 
in the span of $\{v_{2k+1}\, ;\, k\ge0\}$,
 so this quantity is 0. If $d_1=2p\ge 2$ then 
\begin{eqnarray*}
\phi\Big(\Big(1-
\sqrt t \alpha  v_2(s^t_{i_1})\Big)u_{d_1}(s^t_{i_1})\Big)
&=&\phi\Big(\Big(1-
\sqrt t \alpha  v_2(s^t_{i_1})\Big)\Big(
\alpha^p  +\frac 1{\sqrt t}\sum_{j=1}^p \alpha^{p-j}v_{2j}(s_{i_1}^t)\Big)
\Big)\\
&=& \alpha^p - \alpha.\alpha^{p-1}=0
\end{eqnarray*}  
\end{pf}

\begin{cor}
 The map $\rho$ extends to a normal (surjective) 
representation $\Gamma_{t,n}\to 
\Gamma_{1,n}$.
\end{cor}

\begin{pf}
 From the above lemma, it follows that $\tilde \psi$ is positive on 
$C_{t,n}$. As it is normal, it extends to a normal state on $\Gamma_{t,n}$
(that we will denote by $\psi$). Now, the GNS representation of
$\psi$ gives the normal representation.
\end{pf}

\begin{cor}
$\Gamma_{t,n}$ has a direct summand isomorphic to $\Gamma_{1,n}$.
\end{cor}

\begin{rk}\label{compg}
  If the map $\rho$ is not isomorphic, then $\Gamma_{t,n}$ is not a factor. 
Otherwise it is a type $II_1$ factor.
\end{rk}

\begin{rk}
 Denoting by $c^t=(s_1^t)^2+...+(s_n^t)^2$, we have 
$$\eta=\alpha \Big( (n+\frac 1 \alpha)-c^t \Big) \Omega$$
We will study this operator in the next sections.
\end{rk}

\subsection{Case $t<n/(n+\sqrt n)$ and $t>n/(n-\sqrt n)$}

\subsubsection{Case $t<1/2$}

 In this section, we focus on the case $t<1/2$. It is a particular case of the
next one but the arguments are simpler.

When $t<1/2$, the measure of $s^t$ contains one atom at the point $\frac 1
{\sqrt{1-t}} > \sqrt {2t}$. Since in the Calkin algebra
$\rho(s^t)=\rho(\sqrt t s^1)\approx \sqrt t s^1$, we deduce that $C_{t,1}$
contains a compact operator, corresponding to the point projection on
$\frac 1 {\sqrt{1-t}}$ (denoted by $P$) computed on $s^t$. 
Since $\Omega$ is cyclic, this projection is 
one dimensional. From the decomposition $s_i^t \approx s^t
 \oplus (\sqrt t s)_\infty$, we also get that $p=P(s_i^t)$ is a one dimensional
projection. The range of $p$ is the linear span of 
$$\xi_i= \Omega + \frac 1 {\sqrt{1-t}} \sum_{k=1}^\infty \alpha^{\frac {1-k}2}
e_{i^k}
$$
 This vector makes sense as $0<\alpha<1$ when $0<t<1/2$. 

 According to the  Remark \ref{compg} and Proposition \ref{compb},
 on one hand $\Gamma_{t,n}$ ($\infty>n\ge 2$)
 is either $\B(\ell_2(I))$ or not a factor. And  on the other hand it is not 
a factor or is type $II_1$.  

\begin{cor} For $t<1/2$, 
the von Neumann algebra $\Gamma_{t,n}$ is not a factor 
for $2\leq n <\infty$. Moreover it 
has a direct summand isomorphic 
to $\B(\ell_2)$ and another one to $\Gamma_{1,n}$ and
the state $\phi$ is not faithful.
\end{cor}

\begin{lemma}
The vector $\xi_1$ is not cyclic for $\Gamma_{t,n}$.
\end{lemma}
\begin{pf}
It suffices to show that $\eta$ is orthogonal to $\Gamma_{t,n}\xi_1$. But
$\xi_1= P(s_1^t)\Omega$, for $P$ a multiple of the Dirac function 
at $\frac 1{\sqrt {1-t}}$. For any $x\in  \Gamma_{t,n}$, we have
$$ \la x \xi_1, \eta\ra=\la x P(s_1^t)\Omega, \eta\ra=\tau(\rho( x P(s_1^t)))$$
But we know that $\rho(P(s_1^t))=0$ since $P(s_1^t)$ is compact.

\end{pf}

\begin{pf}
The projection $P(s_1^t)$ is of course minimal, so as $\Gamma_{t,n}\xi_1$
is infinite dimensional (for $n\ge 2$), for the $\B(\ell_2(I))$ summand
provided by Proposition \ref{compb}, $I$ is infinite. Then,  
the state $\phi$ can not be faithful because of the  $\B(\ell_2)$ summand.
\end{pf}

\subsubsection{Case $t<n/(n+\sqrt n)$ and $t>n/(n-\sqrt n)$ }

\begin{lemma} For $t\notin \Big[\frac n {n+\sqrt n},\frac n 
{n-\sqrt n}\Big]$, the $C^*$-algebra $C_{t,n}$ contains
a compact operator.
\end{lemma}

{\it\noindent First proof :} It 
is possible to compute explicitly the distribution of 
$c^t=(s^t_1)^2+...+(s_n^t)^2$ for $\phi$
using the $R$-transforms machinery developed in \cite{BLS}. 
We denote by $\gamma_i={s^t_i}^2$, these variables are 
 conditionally free with respect to the distribution of $t{s^1}^2$.
From Theorem 5.2 in  \cite{BLS}, we know the following relations
\begin{eqnarray*}
 R_{c^t}(z)=n R_{\gamma_i}(z) \\
 G_{c^t}(z)=\frac 1 {z-R_{c^t}(G_{tc^1}(z))}
\end{eqnarray*} 
The $R$ and $G$-transforms of $tc^1$ are obtained as usual using freeness
and a change of variable. The computation gives
\begin{eqnarray*}
G_{\gamma_i}(z)=\frac 1 2 \frac{1-2t+\sqrt {1-4t/z}}{z(1-t)-1} \\
R_{\gamma_i}(z)=\frac 1 {1-tz} \\ 
R_{c^t}(z)=\frac n {1-tz} \\
R_{tc^1}(z)=\frac {nt} {1-tz} \\
G_{tc^1}=\frac {(1-n)t+z-\sqrt{(n-1)^2t^2-2(n+1)tz+z^2}}{2tz}\\
G_{c^t}=\frac {(2t-1)z+(1-n)t-\sqrt{(n-1)^2t^2-2(n+1)tz+z^2}}
{2z\Big((t-1)z+n+t(1-n)\Big)}
\end{eqnarray*}  

We can recover the distribution of $c^t$ with respect to $\phi$ from 
the $G$-transform. It turns out that it  has atoms at 
$\frac{n+t(1-n)}{1-t}=n+\frac 1 \alpha$ provided that $t\notin
\Big[\frac n {n+\sqrt n},\frac n {n-\sqrt n}\Big]$.
 
 Actually an eigenvector for the eigenvalue $n+\frac 1 \alpha$ is given by 
$$ \zeta=\sqrt t \Omega + \sum_{k\ge1} \frac 1{(n\alpha)^k} 
\sum_{\substack{|\u i|= 2k\\ i_{2j+1}=i_{2j+2}}} e_{\u i}$$ 

\mk

{\it \noindent Second proof :} We use the Khinchine inequalities
for free products see \cite{H, Bu}. 

Consider the operator 
$$T_k=\sum_{\substack{p\ge 1 \\  k_1+...+k_p=2k\\ k_i \textrm{even} \\ i_1\in 
\{1,...,n\}\\
 i_1\neq i_2\neq..\neq i_p}} u_{k_1}(s^t_{i_1})u_{k_2}(s^t_{i_2})...
 u_{k_p}(s^t_{i_p})$$

In the Calkin algebra, we can identify 
$$\rho(T_k)= \sum_{\substack{p\ge 1 \\ 
 k_1+...+k_p=2k\\ k_i \textrm{even} \\ i_1\in 
\{1,2\}\\ i_1\neq i_2\neq..\neq i_p}} u_{k_1}(\sqrt ts^1_{i_1})u_{k_2}(\sqrt t
s^1_{i_2})...
 u_{k_p}(\sqrt t s^1_{i_p})$$
Thus, we have
$$\rho(T_k)\Omega = \sum_{\substack{|\u i|= 2k\\ i_{2j+1}=i_{2j+2}}} e_{\u i}$$
And $\|\rho(T_k)\Omega\|_2=n^{k/2}$.\\
As $\rho(T_k)\Omega$ is homogeneous of degree $2k$, 
 from the Khinchine inequalities $$\|\rho(T_k)\|\le (2k+1) 
\|\rho(T_k)\Omega\|_2$$ 
But expanding the polynomials using conditional freeness, it comes that  
$\phi(T_k)= n^k \alpha^k$.

So if $n|\alpha|>\sqrt n$ then $\rho$ is not isometric and hence there
is a compact operator in $C_{t,n}$.  
\hfill \qed

\bigskip

\begin{rk}\label{distri} Let 
$$
f(x)=\frac 1 {2x\pi} \frac{\sqrt{\big(x-t(1-\sqrt n)^2\big)
\big(t(1+\sqrt n)^2-x\big)}}
{(t-1)x+n+t(1-n)} \ind_{[t(1-\sqrt n)^2,t(1+\sqrt n)^2]}$$
The distribution of $c^t$ with respect to $\phi$ is 
$$\begin{array}{cl}
f(x) {\rm d}x  & \textrm{ if } t\in \Big[\frac n {n+\sqrt n},\frac n 
{n-\sqrt n}\Big] \\
\displaystyle{f(x) {\rm d}x \,+\,
\frac{(n-1)\big(t-\frac n{n+\sqrt n}\big) \big(t-\frac n{n-\sqrt n} \big)}
{\big(n(1-t)^2+t(1-t)\big)}\, \delta_{n+\frac 1 \alpha}}& \textrm{ if } 
t\notin \Big[\frac n {n+\sqrt n},\frac n 
{n-\sqrt n}\Big]
\end{array}
 $$

\end{rk}
\begin{lemma}
We have $\ker \big( c^t-(n+\frac 1 \alpha)\big)=\C \zeta$.
\end{lemma}

\begin{pf}
We will prove it in several steps. We already know that $\zeta$ is one 
eigenvector for $c^t$. By valuation of a vector in $\mat F_t$,
 we mean  the index of its first non zero 
component according to the natural filtration of $\mat F_t$.
So $\zeta$ has valuation 0. In the following, we let $\xi$ be one eigenvector
(if it exists !) not in $\C \zeta$. We can assume that the valuation of 
$\xi$ is bigger than 1 and that $\xi$ is real.

 From the computation of the $G$-transforms above, the spectrum of 
$tc^1$ is exactly the interval $[t(1-\sqrt n)^2 , t(1+\sqrt n)^2]$ (as 
$\tau$ is faithful).  

\smallskip

{\it\noindent First step : } The valuation of $\xi$ is 1.

\smallskip

 First notice than on tensors of length bigger than 3, $c^t$ acts exactly 
as $tc^1$. Since $\|tc^1\|<n+\frac 1 \alpha$, there is no eigenvector with valuation
bigger than 3.

 Now we focus on the valuation 2. On $\mat H\tens  \mat H$. We consider
the basis given by $f_1=e_{11}+...+e_{nn}, \, f_2=e_{11}-e_{22}, ..., \, 
f_n=e_{11}-e_{nn}$, and the vectors $f_{ij}=e_{ij}$ with $i\neq j$.
 
First the component of $\xi$ on $f_1$ is 0. Indeed, it is clear that 
$\la c^t f_1, \Omega\ra =n$ and for other basis vectors $f$ (and vectors
of valuation bigger than 3), we  have $\la c^t f, 
\Omega\ra =0$. So if $\xi$ has valuation 2, necessarily 
$$ 0=(n+\frac 1 \alpha) \la \xi, \Omega \ra=\la c^t \xi,\Omega\ra=
 \sqrt n \la \xi,f_1\ra$$ 
On the other hand, as above, on vectors of valuation 2 with no component
on $f_1$, $c^t$ acts as $tc^1$. So there is no such vectors.

Consequently, $\dim \ker c^t-(n+\frac 1 \alpha)\le n+1$.
Since $c^t$ is invariant by conjugation with unitaries coming from the 
first quantization. If $\xi$ has valuation 1, we can assume that 
the component of degree $1$ of $\xi$ is $e_1$. 

\medskip

{\it\noindent Second step : } $\xi\in \vect \{ (c^t)^ke_1 \,;\, k\ge 1\}$

\smallskip

 Let $f$ be a continuous function on $\R$ vanishing on the spectrum
of $tc^1$ strictly smaller than 1 except that $f(n+\frac 1 \alpha)=1$.
 As $\rho(c^t)=tc^1$, it follows that $f(c^t)$ is self-adjoint compact.
 Modifying $f$, we can assume that 
$f(c^t)$ is exactly the projection onto $\ker c^t-(n+\frac 1 \alpha)$.  
 
 Then, $f(c^t)e_1$ is non zero as
$$\la f(c^t) e_1,\xi\ra = \la e_1,f(c^t)\xi\ra=  \la e_1,\xi\ra=1$$
So $f(c^t)e_1$ is one eigenvector. It is collinear with $\xi$ 
for its component of degre
1 is collinear to $e_1$ and is non vanishing as there is no eigenvector
of valuation greater than 2.

\medskip

{\it \noindent Third step : } $\xi$ doesn't exist.

\smallskip

It is clear by induction that $\vect \{ (c^t)^ke_1 \,;\, k\ge 0\}=\vect
\{ f_k \; , \; k\ge 0 \}$ with 
$$ f_k= \sum_{\substack{|\u i|= 2k\\ i_{2j+1}=i_{2j+2}}} e_{\u i 1}$$
 We let $\xi=\sum_{k\ge 0} x_k f_k$. 
 we have 
$$\begin{array}{lc}
 c^t f_k = nt f_{k-1} + (n+1)t f_k + f_{k+1} & k\ge 1 \\
 c^t f_0 = (nt+1) f_0 + t f_1&
\end{array}$$
Consequently the sequence $x_i$ has to satisfy the
recursion formula :
$$\begin{array}{lc}
 (n+\frac 1 \alpha) x_0= (nt+1) x_0 + nt x_1 & \\[3pt]
 (n+\frac 1 \alpha) x_k= nt x_{k+1}+(n+1)t x_k + t x_{k-1} & k\ge 1
\end{array}$$
  Hence, $x_k = a (n\alpha)^{-k} + b \alpha^k$, with $b\neq 0$ as soon as
 $n\neq 1$. 
 
But then with this values the series $\sum x_k f_k$
 is not convergent ! So $\xi$ doesn't 
exists
\end{pf}

\begin{thm} For $n\ge 2$ and $t\notin \Big[\frac n {n+\sqrt n},\frac n 
{n-\sqrt n}\Big]$, 
as von Neumann algebras, we have 
$$\Gamma_{t,n}=\mathbb B(\ell_2) \oplus \Gamma_{1,n}$$
Moreover the state $\phi$ is not faithful on $\Gamma_{t,n}$.
\end{thm}

\begin{pf}
 Let $1-q$ be the central support of the representation $\rho$.
As above, we know that $q\neq 0$. We have to show
that $q\Gamma_{t,n}=\mathbb B(\ell_2)$. 

Since $\Omega$ is cyclic for 
$\Gamma_{t,n}$, $q\Omega$ is cyclic in $q\mat F_t$ for $q\Gamma_{t,n}$
(hence non zero).
Let $p$ be the projection onto $\C \zeta$. We have that for $x\in \Gamma_{t,n}$
$$ 0=\tau(\rho(qx))=\la qx\Omega, \eta\ra=\alpha \,\la x\Omega, 
\big(c^t-(n+\frac 1 \alpha)\big)q\Omega\ra$$
So  $(c^t-(n+\frac 1 \alpha))q\Omega=0$, and $q\Omega=\lambda \zeta$ for some 
$\lambda\neq 0$. So $p\le q$, and $q\Gamma_{t,n}$ contains a one dimensional
projection on a cyclic vector so is isomorphic to $\mathbb B(\ell_2)$ 
(as $\Gamma_{t,n} \zeta$ is infinite dimensional). The state can not be 
faithful because of the $\mathbb B(\ell_2)$ summand.
\end{pf}

\subsection{Case $n/(n+\sqrt n)<t<n/(n-\sqrt n)$}

 At the algebraic level, we have already seen that for $i_1\neq...\neq i_l$ 
and $\alpha_j\ge 1$, with $\u i=i_1^{\alpha_1}...i_l^{\alpha_l}$ :
$$ u_{\alpha_1}(s^t_{i_1})...u_{\alpha_{l-1}}(s^t_{i_{l-1}})v_{\alpha_l}(
s^t_{i_l})\Omega= e_{\u i} $$
We define an antilinear map $S$ on $\vect e_{\u i}$ by 
$$ S(e_{\u i})=v_{\alpha_l}(s^t_{i_l})u_{\alpha_{l-1}}(s^t_{i_{l-1}})...
u_{\alpha_1}(s^t_{i_1})\Omega$$

\begin{lemma} The map $S$ satisfies that for any $i_1,...,i_l$
(not necessarily with $i_1\neq...\neq i_l$) and any polynomials $P_1,...,P_l$ :
$$S(P_1(s^t_{i_1})...P_l(s^t_{i_l})\Omega)=\o P_l(s^t_{i_l})... \o 
P_1(s^t_{i_1})\Omega$$
\end{lemma}

\begin{pf}
 Clear by induction on $n$.
\end{pf}

\begin{lemma} For $\frac n {n+\sqrt n}<t<\frac n {n-\sqrt n}$, 
$S$ extends to a bounded operator.
\end{lemma}
\begin{pf}
We  decompose $S$. Let $i_1\neq...\neq i_l$ and $\alpha_j\ge 1$ and define
antilinear operators by :
\begin{eqnarray*}
&A(e_{\u i})=u_{\alpha_l}(s^t_{i_l})u_{\alpha_{l-1}}(s^t_{i_{l-1}})...
u_{\alpha_1}(s^t_{i_1})\Omega& \\
& B(e_{\u i})=u_{\alpha_{l-2}}(s^t_{i_l})u_{\alpha_{l-1}}(s^t_{i_{l-1}})...
u_{\alpha_1}(s^t_{i_1})\Omega & 
\end{eqnarray*}
with the convention that $u_{k}=0$ if $k<0$. Then $S=\sqrt t ( A-\alpha B)$.

From the relations between $v$ and $u$, it follows that $A$ sends tensors of 
length $l$ to a sum of tensors of length $l$, $l-2$, ....
Fix $k\ge 0$, and denote by $A_k$ the component of $A$ that sends tensors of 
length $l$ to tensors of length $l-2k$.  
 Let $J$ be the antiunitary of $\mat F_t$ that reverses 
the order of tensors. Put $C=\sum_{s=1}^n l_1(e_s)^2J$, then   
it is not hard to see that for $|\u i|-|\u j|=2k$ : 
$$\la A_k e_{\u i}, e_{\u j} \ra = \alpha^k f(\u i,\u j) 
\la C^k e_{\u i}, e_{\u j} \ra$$
where $f(\u i,\u j)$ can take the value $1$ or $\frac 1 {\sqrt t}$.
 
Since the coefficients of $C$ are all positive, we get that 
$$\|A_k\|\le \frac{|\alpha|^k}{\sqrt t} \|C^k\|$$
And as $C^*C=n \Id$, we get that 
$$ \|A\|\le \frac 1{\sqrt t} \sum_{k=0}^\infty (\sqrt n |\alpha|)^k$$
which is convergent provided that $\frac n{n+\sqrt n}<t<\frac n {n-\sqrt n}$.

The same kind of arguments  shows that $B$ is bounded.
\end{pf}

Now, we assume that $\frac n{n+\sqrt n}<t<\frac n {n-\sqrt n}$.

\begin{lemma} One has that for any $i_1\neq...\neq i_l$ and $k\le n$
and any polynomials $P_1,...,P_l$ :
$$\Big[SP_1(s^t_{i_1})...P_l(s^t_{i_l})S, s^t_k\Big]=0$$
\end{lemma}
\begin{pf}
 Just a writing game, and the boundedness of $S$
\end{pf}

\begin{cor} For $\frac n{n+\sqrt n}<t<\frac n {n-\sqrt n}$,
 the state $\phi$ is faithful on
$\Gamma_{t,n}$.
\end{cor}

\begin{pf}
$\Omega$ is cyclic for $\Gamma_{t,n}'$ as by the previous lemma 
$S\Gamma_{t,n}S\subset \Gamma_{t,n}'$ and $S$ is invertible (because
$S^2=\Id$).
\end{pf}

\begin{cor}
 For $\frac n{n+\sqrt n}<t<\frac n {n-\sqrt n}$,
 $\Gamma_{t,n}$ does not contain any compact
operator. Moreover the $C^*$-algebras $C_{t,n}$ and $C_{1,n}$ are isomorphic.
\end{cor}

\begin{pf}
 If $\Gamma_{t,n}$ contains a compact operator, then  
as above there is a  direct summand of $\Gamma_{t,n}$ isomorphic to 
$\B(\ell_2(I))$. But as $\Omega$ is separating, $I$ must consist of one 
point, which is impossible as $s^t_1$ as no eigenvector.

The second assertion follows from the first one since the Calkin map
is then an isomorphism in both case.  
\end{pf}

\begin{thm} For  $\frac n{n+\sqrt n}<t<\frac n {n-\sqrt n}$,
the von Neumann algebra $\Gamma_{t,n}$ is isomorphic to $\Gamma_{1,n}$. 
\end{thm}

\begin{pf}
 One just needs to show that the map $\rho$ is faithful. To do so it suffices
to show that $\psi$ is faithful.
 
Let $b=S(1-\sqrt{t} \alpha \sum_{i=1}^n v_2(s_i^t))S\in \Gamma_{t,n}'$.
We have for $x\in \Gamma_{t,n}$
$$\psi(x)=\la b^*x \Omega,\Omega\ra.$$
Applying it to $x=y^*y$ and as $\Omega$ is cyclic it follows that 
$b$ is positive. Say $b=a^2$, then 
$$\psi(x^*x)=\|xa \Omega\|^2$$
 Moreover the distribution of $b$ for $\phi$ is absolutely
continuous with respect to the Lebesgue measure (see Remark \ref{distri}), 
it follows that there is a net 
of elements  $(c_i)$ in the von Neumann generated by $b$ so that
$c_ia\Omega\to \Omega$. Now, if for $x\in \Gamma_{t,n}$, $\psi(x^*x)=0$ then
$$0=\lim c_ixa\Omega= \lim xc_ia\Omega=x\Omega.$$
So $x=0$ as $\Omega$ is cyclic.
\end{pf}

\begin{rk}
 The Tomita-Takesaki operator $S$ is bounded on $\Gamma_{t,n}$.
\end{rk}

\begin{rk}\label{inversion}
It is possible to reverse all the arguments above. To do so, define
a normal linear form $\tilde \phi$ on $\Gamma_{1,n}$ by
$$\tilde \phi(x)=\la x \Omega, \sum_{k\ge 0}\alpha^k 
\sum_{\substack{|\u i |=2k\\ i_{2j+1}=i_{2j+2}}} e_{\u i} \ra$$
It is well defined because the vector on the right side makes sense 
as $\sqrt n |\alpha|<1$ \\
 It is clear that on the set of polynomials in $s^1_i$, it coincides with
$\phi$. By continuity, we get that on $\Gamma_{t,n}$, we have
$$\phi = \tilde \phi \rho$$
   Then the GNS construction of $\Gamma_{1,n}$ for $\tilde \phi$ gives 
$\rho^{-1}$.
\end{rk}

\subsection{Case $t=n/(n\pm \sqrt n)$  }\label{caseg}

 In this situation, we will adopt the strategy of the last remark.

 Let $\mathcal P$ be the set of noncommutative polynomials in the letters 
$X_1,...,X_n$. It has a natural $*$-algebra structure. \\
  For $0<r<1$, we define a linear functional $\phi_r$ on $\mathcal P$ : for 
$i_1\neq...\neq i_l$ and $\alpha_j\ge 1$ 
$$ \phi_r\Big(u_{\alpha_1}(X_{i_1})...u_{\alpha_l}(X_{i_l})\Big)=\left\{
\begin{array}{ll}
(r^2 \alpha)^{\sum \alpha_i/2} &  \textrm{ if all $\alpha_i$ are even} \\
0  &  \textrm{ otherwise}
\end{array}\right.$$
 Formally, we can identify $\mathcal P$ with the $*$-algebra generated
by $\sqrt t s^1_i$ in $C_{1,n}$ (say $\pi(X_i)=\sqrt t s^1_i$) and to 
$*$-algebra generated
by $ s^t_i$ in $C_{t,n}$ (say $\sigma(X_1)= s^t_i$). \\
 Let $T_r$ be the second quantization associated to $r\Id$ on $\Gamma_{1,n}$,
it is unital completely positive and let $\pi(\mathcal P)$ invariant. \\
The functional $\phi_r$ is made so that 
$$\phi_r(P)=\phi\Big( \sigma(\pi^{-1}(T_r(\pi(P))))\Big)$$ 
Indeed, for $i_1\neq...\neq i_l$ and $\alpha_j\ge 1$ 
\begin{eqnarray*}
\pi^{-1}(T_r(\pi(u_{\alpha_1}(X_{i_1})...u_{\alpha_l}(X_{i_l}))))&=&
\pi^{-1}(T_r(u_{\alpha_1}(\sqrt ts^1_{i_1})...u_{\alpha_l}(\sqrt ts^1_{i_l})))
\\ &=& \pi^{-1}(
r^{\sum \alpha_i} u_{\alpha_1}(\sqrt ts^1_{i_1})...u_{\alpha_l}
(\sqrt ts^1_{i_l})) \\ &=&
r^{\sum \alpha_i}u_{\alpha_1}(X_{i_1})...u_{\alpha_l}(X_{i_l})
\end{eqnarray*}
and by conditional freeness :
$$\phi\Big( \sigma(u_{\alpha_1}(X_{i_1})...u_{\alpha_l}(X_{i_l}))\Big)=
\phi\Big(u_{\alpha_1}(s^t_{i_1})...u_{\alpha_l}(s^t_{i_l})\Big)=
\left\{
\begin{array}{ll}
 \alpha^{\sum \alpha_i /2} &  \textrm{ if all $\alpha_i$ are even} \\
0  &  \textrm{ otherwise}
\end{array}\right.$$

\begin{lemma} The functional $\phi_r$ is positive on $\mathcal P$.
 \end{lemma}

\begin{pf}
 The $*$-algebra $A_i$ generated by $X_i$ in $\mathcal P$ has two particular
functionals on it. The first one is the restriction of $\phi_r$ (still denoted
by $\phi_r$) and the second one is $\psi\pi$. From the definition of $\phi_r$
the algebra $(A_i,\phi_r,\psi\pi)$ are conditionally free. Hence the result
will follow from \cite{BLS}
 Theorem 2.2 provided that we can show the positivity 
of $\phi_r$ and $\psi\pi$ on each $A_i$. For $\psi\pi$ this is obvious as
$\psi$ is a state. For $\phi_r$, it suffices to notice that the distribution
of $X_i$ for $\phi_r$ the distribution of $rs^t_i$ for $\phi$, so comes from
a probability measure. 
\end{pf}

\begin{cor} The functional $\phi_r\pi^{-1}$ extends to a normal state
on $\Gamma_{1,n}$ (denoted by $\phi_r$).
\end{cor}

\begin{pf}
 By the preceding Lemma, it is positive on $\pi(\mathcal P)$ which is 
weak-$*$ dense in $\Gamma_{1,n}$. The result follows from the representation
$$ \phi_r(\pi^{-1}(x))=\la x \Omega, \sum_{k\ge 0}(r^2\alpha)^k 
\sum_{\substack{|\u i |=2k\\ i_{2j+1}=i_{2j+2}}} e_{\u i} \ra$$
It is well defined because the vector on the right side makes sense 
as $r^2 \sqrt n |\alpha|=r^2<1$.
\end{pf}

\begin{thm}
 The $C^*$-algebras $C_{n/(n\pm\sqrt n),n}$ 
and $C_{1,n}$ are isomorphic and the state $\phi$ is faithful
on $C_{n/(n\pm\sqrt n),n}$. 
\end{thm}

\begin{pf}
 There is a state on $C_{1,n}$ defined by $\tilde \phi=
\lim_{r\to 1,\mathfrak U} \phi_r$, where the limit is taken along some 
non trivial ultrafilter. It is 
clear from the formulas  on $\mathcal P$, that for 
$x\in \C_{t,n}$, we have 
$$ \phi(x)= \tilde \phi ( \rho(x)).$$

 We consider the GNS construction  of $\tilde \phi$ for $C_{1,n}$. The 
underlining Hilbert space is exactly $\mat F_{n/n\pm\sqrt n}$ from the above
identities, so this GNS construction is the inverse of the Calkin map $\rho$.

 Assume $x\in C_{n/(n\pm\sqrt n),n}^+$ is such that $\phi(x^2)=0$, then
$x\Omega=0$. It follows that $\psi(\rho(x))=\la x\Omega,\eta\ra=0$. So 
$\rho(x)=0$ as $\psi$ is faithful on $\Gamma_{1,n}$, so $x=0$. 
\end{pf}

 The analogue result for the von Neumann algebras is a consequence of the 
next section. In particular they show that the map $\tilde \phi$ is actually 
normal.

\subsection{Case $n=\infty$}
 
 This case was already treated in \cite{W} 
and was a starting point to this work.

\smallskip

If $t\neq 1$ then it can be checked that 
$$\frac 1 k \sum_{i=1}^k (s^t_i)^2 
\mathop{\longrightarrow}\limits_{k\to \infty}
(1-t) P + t  \Id$$
in the weak-* topology, 
where $P$ is the orthogonal projection onto $\Omega$. As $\Omega$ is cyclic,
Proposition \ref{compb} ensures that 
$$\Gamma_{t,\infty}=\B(\mat F_t)$$

\section{Generalizations}

 In this section, we explain how to extend the previous results.

 Let $I_i$ be bounded measurable subsets of $\R$, 
$A_i=L_\infty(I_i,\mu_i)$ be $n$ commutative 
von Neumann algebras. The integration with respect to $\mu_i$ will be called
$\phi_i$ and $s_i$ is the identity function  on $I_i$.
 Assume that we are given $\psi_i$,   distinguished normal states
on $A_i$. As a consequence, $\psi_i$ has a density 
$f_i$ with respect to $\phi_i$, we put $\b f_i=f_i-1$. 
We will denote by $v_k^i$ and $u_k^i$ the 
orthogonal polynomials for $\phi_i$ and $\psi_i$.

 We consider the algebraic conditional free product $\tilde 
A=*_{i=1}^n (A_i,\phi_i,\psi_i)$.
The state $\phi$ is the conditional free product state and $\psi$ is the 
free product of the $\psi_i$.
Then we can talk about the von Neumann reduced conditional free product, 
corresponding to GNS construction of $\phi$, we denote it by $A$
$$A \subset \B(\mat F) \qquad \phi=\la . \Omega,\Omega\ra $$
 $\mat F$ is the conditional Fock space and $\Omega$ the vacuum state, 
we refer to 
 \cite{Ml} Section 6 for more details about this construction. One of the main
trouble comes from the fact 
 that in general the state $\phi$ is not faithful on $A$
(for $t$-gaussians for instance).
Nevertheless, the injections $i_i:A_i\to A$ are normal and isometric and 
preserve the states (i.e. $i_i^*(\phi)=\phi_i$, so we will drop the 
indexes) (this corresponds to Proposition \ref{basic}). Unfortunately, it 
seems that in general there is no conditional expectation from $A$ to $A_i$
which is state preserving.  

 In this situation, the Lemmas \ref{ident} and \ref{densite} remain true.
For $i_1\neq...\neq i_l$ 
and $\alpha_j\ge 1$, let $\u i=i_1^{\alpha_1}...i_l^{\alpha_l}$ and define :
$$ e_{\u i}=u_{\alpha_1}^{i_1}(s_{i_1})...u_{\alpha_{l-1}}^{i_{l-1}}
(s_{i_{l-1}})v_{\alpha_l}^{i_l}(s_{i_l})\Omega$$
\begin{lemma}\label{baseorth}
 The family $(e_{\u i})_{\u i\in \{1,...,n\}^{<\infty}}$
 is an orthonormal basis of $\mat F$. \\
 If the densities $f_i$ are bounded (or in $L^2(\mu_i)$), the state $\psi$ on 
the algebraic free product extends to a normal state
of $A$.
\end{lemma}

\begin{pf}
To prove the first assertion we will need to evaluate expressions of the form
$\phi(a b_1... b_l d )$, where $a\in A_{i_0}$, $b_k\in A_{i_k}$, 
$\psi_{i_k}(b_k)=0$ and 
$d\in A_{i_{l+1}}$. With $b=b_1... b_l$ and $x=\phi_{i_1}(b_1)...
\phi_{i_l}(b_l)$ 

\begin{eqnarray*}
 \phi( a b d)&=&\phi \big( (a -\psi_{i_0}(a))b(d -\psi_{i_{l+1}}(d))\big) +
\psi_{i_0}(a)\phi\big( b(d-\psi_{i_{l+1}}(d))\big)+\\ & & 
\phi \big( (a -\psi_{i_0}
(a))b\big)
\psi_{i_{l+1}}(d)+ 
\psi_{i_0}(a)\phi(b)\psi_{i_{l+1}}(d)\\ &=& \big(\phi(a)-\psi_{i_0}(a)\big)x
\big(\phi(d)-\psi_{i_{l+1}}(d)\big) \psi(a) x \big(\phi(d)-\psi_{i_{l+1}}(d)
\big)+\\ & & +
\big(\phi(a)-\psi_{i_0}(a)\big)x\psi_{i_{l+1}}(d) +
 \psi_{i_0}(a)x\psi_{i_{l+1}}(d) \\ &=& \phi(a)x
\phi(d)\\ &=& \phi(a)\phi(b)\phi(d)
\end{eqnarray*}   

For $\u i=i_1^{\alpha_1}...i_l^{\alpha_l}$ and 
$\u j=j_1^{\beta_1}...j_m^{\beta_m}$, we have 
$$\la e_{\u i} , e_{\u j} \ra = \phi\Big(
v_{\beta_m}^{j_m}(s_{j_m})
u_{\beta_{m-1}}^{j_{m_{l-1}}}(s_{j_{m-1}})...u_{\beta_1}^{j_1}
(s_{j_{1}})
u_{\alpha_1}^{i_1}(s_{i_1})...u_{\alpha_{l-1}}^{i_{l-1}}
(s_{i_{l-1}})v_{\alpha_l}^{i_l}(s_{i_l})\Big)$$
 So if $i_1\neq j_1$ or $\alpha_1\neq\beta_1$ then  as 
$\phi(v_{\alpha_l}^{i_l}(s_{i_l}))=0$, the previous identity ensures that this
scalar product is 0. In the remaining case $i_1=j_1$ and $\alpha_1=\beta_1$,
so $\psi_{i_1}(u_{\beta_1}^{j_1}
(s_{j_{1}})
u_{\alpha_1}^{i_1}(s_{i_1}))=1$ and the above identity gives
$$\la e_{\u i} , e_{\u j} \ra=\la e_{i_2^{\alpha_2}...i_l^{\alpha_l}} , 
e_{j_2^{\beta_2}...i_l^{\beta_m}} \ra$$
and an induction completes the proof as the $v_n^i$ are orthonormal.

  If the $f_i$ are bounded, then 
$$c=1+ \sum_{i=1}^n\b f_i(s_i)\in A$$
and for all $x\in \tilde A$,
$$\phi(cx)=\psi(x)$$
 By linearity and freeness, it suffices to check it for 
$x=u_{\alpha_1}^{i_1}(s_{i_1})...u_{\alpha_{l-1}}^{i_{l-1}}
(s_{i_{l-1}})u_{\alpha_l}^{i_l}(s_{i_l})$  with $l\ge 1$, $\alpha_k>0$ and
$i_1\neq...\neq i_l$ as $\phi(c)=1$. 
We have
\begin{eqnarray*}
\phi(cx)& = & \phi\Big(\big(1+ \sum_{i=1}^n\b f_i(s_i)\big)
u_{\alpha_1}^{i_1}(s_{i_1})...u_{\alpha_{l-1}}^{i_{l-1}}
(s_{i_{l-1}})u_{\alpha_l}^{i_l}(s_{i_l})\Big)\\
& = & \phi\Big( f_{i_1}(s_{i_1})u_{\alpha_1}^{i_1}(s_{i_1})
...u_{\alpha_{l-1}}^{i_{l-1}}
(s_{i_{l-1}})u_{\alpha_l}^{i_l}(s_{i_l})\Big) \\ & = &
\phi\big( f_{i_1}(s_{i_1})u_{\alpha_1}^{i_1}(s_{i_1})\big) \phi(
u_{\alpha_2}^{i_2}(s_{i_2}))...\phi(u_{\alpha_l}^{i_l}(s_{i_l}))\\
&=& \psi( u_{\alpha_1}^{i_1}(s_{i_1})) \phi(
u_{\alpha_2}^{i_2}(s_{i_2}))...\phi(u_{\alpha_l}^{i_l}(s_{i_l}))\\
&=& 0
\end{eqnarray*}
This gives the normality and the extension of $\psi$ as $c\in A$. 

If the $f_i$ are only in $L^2$, then we can consider
$$\eta=\Big(1+ \sum_{i=1}^n\b f_i(s_i)\Big).\Omega$$
and simply replace $\phi(cx)$ by $\la x\Omega,\eta\ra$ in the above
proof.
\end{pf}

\begin{rk}
 We need the assumption on the densities because of the lack of 
conditional expectation. So we can not define $c$ correctly when 
$f_i$ are only in $L^1$. It is very likely that in concrete situations
one can find an appropriate proof of this result.
\end{rk}

\begin{cor} If $f_i\in L^2(\mu_i)$, 
the conditional free product $*_{i=1}^n (A_i,\phi_i,\psi_i)$
 has a direct summand isomorphic to the free product 
$*_{i=1}^n (A_i,\psi_i)$.
\end{cor}

\begin{rk}
In order to have isometric copies of $A_i$ in the free product
$*_{i=1}^n (A_i,\psi_i)$, one needs to assume that the GNS construction 
of $A_i$ for $\psi_i$ to be faithful. So it boils down to assume that 
$\phi_i$ is absolutely continuous with respect to $\psi_i$, this corresponds
to the case $t\geq \frac 1 2$ for $t$-gaussians.
\end{rk}

\begin{rk}
 The basis $e_{\u i}$ is natural in the sense that on it, the generators 
acts like gaussian operators (creation plus annihilation). In general, one 
does not recover interacting Fock spaces.
\end{rk}

\begin{thm}
Assume that the densities $f_i$ are bounded and that the distribution of
$c$ with respect to 
 $\phi$ does not have atom at 0. Then as von Neumann algebras
$$*_{i=1}^n (A_i,\phi_i,\psi_i)=*_{i=1}^n (A_i,\psi_i)$$
If moreover $\psi_i$ is faithful on $A_i$ for $i=1,...,n$ then
the state $\phi$ is faithful.
\end{thm}

\begin{pf}
 We denote by $\rho$ the representation from the conditional free product
to the free product.

 The GNS representation for 
$\phi$ restricted to the von Neumann generated by $c$ is included in 
the reduced conditional free product. Then the  assumption on $c$ ensures that
one can find $b_n$ (self-adjoint) 
in the $C^*$-algebra generated by $c$ so that 
$b_nc\Omega \to \Omega$. 

Let $x$ be such that $\rho(x)=0$.
  Hence, for any $a$ and $d$ in $A$, 
$$0=\lim_{k\to \infty} \psi(\rho(b_naxd))= \lim_{k\to \infty} \phi(cb_naxd)
=\lim_{k\to \infty}\la axd \Omega, b_nc\Omega\ra=\la xd \Omega, a^*\Omega\ra=
 \phi(axd)$$
 So $x=0$ as $\Omega$ is cyclic for $A$.

 If the $\psi_i$ are faithful, then their free product $\psi$ 
is also faithful by a result of Dykema \cite{D}. 
Assume $x$ is so that $\phi(x^*x)=0$, then 
by the Cauchy-Schwarz inequality, $\phi(cx^*x)=0$. So 
one has $\psi(\rho(x^*x))=0$ and $\rho(x)=0$. Consequently
 $x=0$ as $\rho$ is one to one.
\end{pf}
 
\begin{rk} If the distribution of $c$ has atom at 0,
 then the reduced conditional
free product is not a factor.
\end{rk}

 Another approach, to show that the free product and the conditional free 
products coincide, is to adopt the strategy of the Remark \ref{inversion}
or Section \ref{caseg}. Unfortunately, it seems that it is not as good as 
the above Theorem, as one can not use it to get the limit cases 
$t=\frac n{n\pm \sqrt n}$. 
 
 We let $f_i$ be the natural orthonormal basis in the free Fock space
associated to the above distributions. Consider the vector
$$\zeta= \sum_{\u i=i_1^{\alpha_1}...i_l^{\alpha_l}} \phi(u_{\alpha_1}(s_{i_1})
...u_{\alpha_l}(s_{i_l})) f_{\u i}$$
 It exists if and only if 
$$  \sum_{\u i=i_1^{\alpha_1}...i_l^{\alpha_l}} \phi_{i_1}
(u_{\alpha_1}(s_{i_1}))^2...\phi_{i_l} (u_{\alpha_l}(s_{i_l}))^2<\infty.$$
 Then,  the normal state on $*(A_i,\psi_i)$ defined by 
$$\tilde \phi (x) =\la x\Omega,\zeta\ra$$
coincides with $\phi$ on $\mat P$. So the GNS construction of this normal 
state gives a representation  from  the free product to the conditional free
product.  

\medskip 

{\it\noindent Acknowledgement :} The author would like to thank 
 the Wroc{\l}aw department of mathematics 
for introducing him to the subject
 and also for their kind hospitality.

\end{document}